\documentclass{amsart} 

\usepackage{amsmath} 
\usepackage{amssymb} 

\usepackage{amsfonts}

 \usepackage{tikz}
\usepackage{pgfplots}

\begin{document}

\newcommand{\BB}{{\mathbb B}}
\newcommand{\CC}{{\mathbb C}}
\newcommand{\GG}{{\mathbb G}}
\newcommand{\HH}{{\mathbb H}}
\newcommand{\PP}{{\mathbb P}}
\newcommand{\QQ}{{\mathbb Q}}
\newcommand{\RR}{{\mathbb R}}
\newcommand{\TT}{{\mathbb T}}
\newcommand{\ZZ}{{\mathbb Z}}

\newcommand{\cA}{{\mathcal A}}
\newcommand{\cC}{{\mathcal C}}
\newcommand{\Coh}{{\mathcal Coh}}
\newcommand{\cD}{{\mathcal D}}
\newcommand{\cE}{{\mathcal E}}
\newcommand{\cF}{{\mathcal F}}               
\newcommand{\cI}{{\mathcal I}}
\newcommand{\cH}{{\mathcal H}}
\newcommand{\cM}{{\mathcal M}}
\newcommand{\cO}{{\mathcal O}}
\newcommand{\cS}{{\mathcal S}}
\newcommand{\cSets}{{\mathcal Sets}}
\newcommand{\cSub}{{\mathcal Sub}}
\newcommand{\cT}{{\mathcal T}}
\newcommand{\cZ}{{\mathcal Z}}

\newcommand{\ch}{\mbox{ch}}
\newcommand{\coker}{\mbox{coker}}
\newcommand{\Ext}{\mbox{Ext}}
\newcommand{\Hilb}{\mbox{Hilb}}
\newcommand{\Hom}{\mbox{Hom}} 
\newcommand{\rH}{\mbox{H}} 
\newcommand{\im}{\mbox{Im}}
\newcommand{\length}{\mbox{length}}
\newcommand{\NS}{\mbox{NS}}
\newcommand{\ord}{{\rm ord}}
\newcommand{\rk}{\mbox{rk}}
\newcommand{\re}{{\rm Re}}
\newcommand{\Quot}{{\rm Quot}}
\newcommand{\Sym}{\mbox{Sym}}

\newcommand{\id}{\mbox{id}}
\newcommand{\tr}{\mbox{tr}}
\newcommand{\End}{\mbox{End}}
\newcommand{\Aut}{\mbox{Aut}}
\newcommand{\Perm}{{\rm Perm}}
\newcommand{\sgn}{{\rm sgn}}
\newcommand{\PGL}{{\rm PGL}}

\newcommand{\Sec}{{\rm Sec}}

\newcommand{\fun}{\rightarrow}
\newcommand{\nt}{\noindent}

\begin{center}

{\bf Secants, Socles and Stability}

\bigskip

Aaron Bertram\footnote{Partially supported by the Simons Foundation through a Travel Grant} and Brooke Ullery

\medskip

{\it Dedicated to Robert Lazarsfeld} 

\end{center}

\medskip

\section{Introduction} 

Let $k$ be an algebraically closed field of characteristic zero and:
\[ S_\bullet = k[x_0,...,x_n]_\bullet \] 
be the polynomial ring, graded by degree.  The geometry of the secants to:
\[ \PP(S_1) =: \PP^n \hookrightarrow \PP^n_d := \PP(S_d) \] 
(the ``$d$-uple embeddings'') is a polynomial version of the Waring problem for natural numbers \cite{War}.
In this context, one is interested in the secant loci in $\PP^n_d$ swept out 
by $m$-planes spanned by $m+1$ points of $\PP^n$ (possibly including  limits), and one wants to know which points in such a secant locus 
``remember'' the $m+1$ points. 

\medskip

Alternatively, one can view $\PP(S_d)$ as the parameter space for degree $d$ socles $\sigma$ of Gorenstein 
graded $S$-modules $R_\bullet$. From that point of view, there are a pair of important numerical invariants
associated to a point $\sigma \in \PP(S_d)$, namely, 

\medskip

(i) the palindromic {\bf Hilbert function} of $R_\bullet$, and 

\medskip

(ii) the (also palindromic) {\bf betti table} of the minimal free resolution of $R_\bullet$. 

\medskip

The second refines the first, and  the minimal free resolution is also filtered by ``linear strands.''
Some guidance is needed, however, in finding the right coarsening of this overly fine filtration to explain, for example, 
why the classical variety of secant lines to the Veronese (two-uple) embedding of $\PP^n$ is defective (of smaller than the expected dimension) while the varieties
of secant lines to the three-uple (and higher) embeddings are non-defective. 

\medskip

In this paper, we appeal for such guidance to the canonical Harder-Narasimhan filtrations 
with respect to two natural {\it stability conditions} on the bounded derived category ${\mathcal D}^b(\PP^n)$ of coherent sheaves
on $\PP^n$. This involves a third interpretation of the projective space $\PP^n_d$ via the vector space: 
\[ \sigma \in \Hom_{\PP^n}({\mathcal O}_{\PP^n}, \omega_{\PP^n}(-d)[n]) \] 
in the derived category, recovering the interior of the minimal free resolution of $R_\bullet$ via the {\bf cone}  over the morphism. From
this point of view, the uniqueness of the ``spanner'' of a secant reduces to the uniqueness of Harder-Narasimhan filtrations, 
and the challenge is to see how good the stability condition is at distinguishing betti tables from one another. 

\medskip

This approach generalizes, in principle, to other varieties embedded by complete linear series. We will not discuss this here, but our treatment for $\PP^n$ mirrors the use of  
classical {\it Mumford (in)stability} to distinguish secant varieties to embedded curves, as explained to each of the authors by their advisor, Robert Lazarsfeld. \cite{Be1,Be2,Ull}

\section{Gorenstein Quotient Rings.} 

\nt {\bf Definition 2.1.} A finite-dimensional graded quotient ring:
\[ S_\bullet \rightarrow R_\bullet = S_\bullet/I_\bullet \] 
is {\bf Gorenstein} with socle in degree $d$ if: 

\medskip

(a) $\dim_k R_d = 1$

\medskip

(b) for all $e\in \ZZ$, the multiplication pairing: 
\[ R_e \times R_{d - e} \rightarrow R_d \] 
is perfect, inducing vector space isomorphisms $R_{d-e} \cong R_e^\vee$. 

\medskip

\nt {\bf Note.} The ``socle'' $I_d = \ker(S_d \rightarrow R_d)  \in (S_d)^\vee$ determines $I_e$ for all $e$ via:
\[ I_e = \ker(S_e \rightarrow S_{d-e}^\vee ) \ \mbox{for the map} \ (F,G) \mapsto F\cdot G + I_d \in R_d    \]  
and $I_e = S_e$ if and only if $e$ is outside the range $[0,d]$. 

\medskip

\nt {\bf Example 2.1.} For socles $\sigma = I_d$ for which $\dim(R_1) = 1$, 
\[ I_\bullet = \langle I_1 \cdot S_\bullet, G \rangle \ \mbox{for any} \ G \in S_{d+1} - I_1 \cdot S_d \]  
and $\dim(R_\bullet) = (1,1,....,1)$. These socles therefore correspond to 
$I_1  \in (S_1)^\vee$. 

\medskip

Passing to the dual, we may in general let: 
\[ V = (S_1)^\vee \] 
so that $\PP^n$ is the projective space of {\it lines} in $V$, and then
\[ (S_d)^\vee \ \mbox{is the space of symmmetric $d$-tensors} \ \Sym^dV \subset V^{\otimes d} \]  
and the $d$-uple embedding is the map $\PP^n \rightarrow \PP^n_d \ \mbox{given by} \ v \mapsto v^d$. Moreover:

\medskip

$\bullet$ The image of the $d$-uple embedding is the locus of socles with $\dim_k(R_1) = 1$. 

\medskip

The Waring problem asks when a symmetric $d$ tensor $t$ has the form: 
\[ t = \sum _0^m c_i v_i^d \] 
for distinct vectors $v_i \in V$. This is exactly the problem of determining whether a socle $\sigma = t$
is in the locus of secant spaces to $m+1$ points of the $d$-uple embedding. 

\medskip

\nt {\bf Example 2.2.} When $d = 2$, there is an identification: 
\[ (S_2)^\vee = \{ q: S_1 = V^\vee \rightarrow V = (S_1)^\vee \ | \ q = q^T\} \] 
with the symmetric matrices, and the indeterminate middle dimension:
\[ \dim(R_\bullet) = (1, \ r, \ 1) \] 
is the {\bf rank} of the matrix $q$ associated to $\sigma$ since
$I_1 = \ker(q)$.

\medskip

The Waring problem in this context translates into the fact that a symmetric matrix of rank 
$r$ is a always a sum of $r$ distinct symmetric matrices of rank one, but of course the 
choice of the $r$ matrices is not unique as, for example, 
\[ v_0^2 + v_1^2 = \frac 12 (v_0 + v_1)^2 + \frac 12 (v_0 - v_1)^2 \] 

\nt {\bf Example 2.3.} When $n = 1$ and $d$ is arbitrary, the Waring problem has a very satisfactory answer. In that case, 
 the Hilbert function of $R_\bullet$ is a trapezoid: 
\[ \dim(R_e) = \left\{ \begin{array}{l} e + 1 \ \mbox{for} \ 0 \le e \le m - 1 \\  \\ 
m  \ \mbox{for} \ m - 1 \le e  \le d- (m - 1) \\ \\
d - e + 1 \ \mbox{for} \ d- (m-1)  \le e \le d  \end{array} \right.\] 
and Waring's problem (modified to allows for multiplicites) has the nicest possible result
when $2m < d$. Namely,  
each linear combination $t = \sum_{i  = 0}^m  c_iv_i^d$ of distinct rank one tensors 
uniquely determines the vectors $v_0,....,v_m$ (up to scalar multiples). 

\medskip

Thus, for example, when $n = 1$ and $d = 3$ and $m = 1$, each sum: 
\[ v_0^3 + v_1^3  \] 
has no other expression as a sum of two cubes. 

\medskip

\nt {\bf Note.} The span of a non-reduced subscheme $Z \subset \PP^1$ of length $m+1$ 
may be absorbed into the Waring problem. When $d = 3$ and $m = 1$, the span of double points is the union of the tangent 
lines to the twisted cubic in $\PP^3$, and is disjoint from the variety of secant lines to pairs of points (minus the points themselves). 
The ``right'' way to describe the $n = 1$ case is to say that when $2m \le d $,

\medskip

(1) each linear span of a subscheme $Z \subset \PP^1$ of length $m+1$  is a $\PP^{m} \subset \PP^{d}$. 

\medskip

(2) the strata $X_m - X_{m-1}$ of spans of  length $m+1$ subschemes are locally closed.

\medskip

(3) when $2m < d$, each $\sigma \in X_e - X_{e-1}$ determines the subscheme that spans it. 

\medskip

Next, consider the {\bf minimal free resolution} of a Gorenstein quotient:
\[ \cdots \rightarrow F_2 \rightarrow F_1 \rightarrow S_\bullet \rightarrow R_\bullet \rightarrow 0 \] 
where $F_1$ is the free module minimally generating $I_\bullet$, 
$F_2$ is free on a minimal set of relations, etc. By the Hilbert Syzygy Theorem, this terminates at $F_{n+1}$. 

\medskip

\nt {\bf Example 2.4.} When $d = 0$, the ideal is maximal:
\[ 0 \rightarrow \langle x_0,...,x_n \rangle \rightarrow S_\bullet \rightarrow k \rightarrow 0 \] 
and the minimal free resolution is the {\bf Koszul complex}:
\[ 0 \rightarrow \wedge^{n+1} V \otimes S(-n-1) \rightarrow \cdots \rightarrow \wedge^iV \otimes S(-i) \rightarrow \cdots \rightarrow V \otimes S(-1) \rightarrow S \rightarrow k \] 

\nt {\bf Definition 2.2 .} The array of natural numbers $b_{i,j}$ such that:
\[ F_i = \bigoplus S(-j)^{b_{i,j}}  \] 
is the betti table of $R_\bullet$, usually arranged (since $j \ge i$) as: 

\[ \begin{tabular}{|c|c|c|c|c|c|} 
\hline $b_{0,0} = 1$ & $b_{1,1}$ & $b_{2,2}$  & $b_{3,3}$ & $\cdots$ & $b_{n+1,n+1}$ \\ \hline 
& $b_{1,2}$ & $b_{2,3}$ & $b_{3,4}$ & $\cdots$ & $b_{n+1,n+2}$  \\ \hline 
& $\vdots$ & $\vdots$ & $\vdots$ & $\vdots$ & $\vdots$ \\ \hline 

& $b_{1,d+1}$ & $b_{2,d+2}$ & $b_{3,d+3}$ & $\cdots$ & $b_{n+1,n+1+d}$ \\ \hline 

\end{tabular} \] 
 
 \medskip
 
 \nt so that, for example, the betti table of the Koszul complex has a single row:

\[ \begin{tabular}{|c|c|c|c|c|c|} 
\hline $b_{0,0} = {{n+1} \choose 0}$ & ${{n+1} \choose 1}$ & ${{n+1} \choose 2}$  & ${{n+1} \choose 3} $ & $\cdots$ & $b_{n+1,n+1} = {{n+1} \choose {n+1}}$  \\ \hline 

\end{tabular} \] 

\medskip

The left-right symmetry of this Koszul row generalizes via a well-known: 

\medskip

\nt {\bf Duality.} \cite{Mac}  If $R_\bullet$ is Gorenstein with socle in degree $d$, then 

\medskip

(i) $F_{n+1} = S(-n-1-d)$, i.e. $b_{n+1,n+1+d} = 1$ and $b_{n+1,n+1+e} = 0$ for $e < d$.   

\medskip

(ii) $F_{n+1-m} \cong \mbox{Hom}_S(F_m,S)(-n-1-d)$, and when applied to the complex, 

\medskip

(iii) $\mbox{Hom}_S(F_*,S)(-n-1-d) = -F_{n+1-*}$  (skew symmetry)

\medskip

\nt {\bf Corollary.} 
$ b_{n+1,n+1} = \cdots = b_{n+1,n +1 + d - 1} = 0$. 

\medskip

\nt {\bf Example 2.5.} When $d = 1$, there is a (smaller) Koszul complex for the first row: 
\[ \begin{tabular}{|c|c|c|c|c|c|c|} 
\hline $b_{0,0} = {{n} \choose 0}$ & ${{n} \choose 1}$ & ${{n} \choose 2}$  & ${{n} \choose 3} $ & $\cdots$ & $b_{n,n} = {{n} \choose {n}}$ & $0$ \\ \hline 

\end{tabular} \] 

\medskip

\nt to which the second row is appended according to the Theorem: 
\[ \begin{tabular}{|c|c|c|c|c|c|c|} 
\hline $b_{0,0} = {{n} \choose 0}$ & ${{n} \choose 1}$ & ${{n} \choose 2}$  & ${{n} \choose 3} $ & $\cdots$ & ${{n} \choose {n}}$ & $0$  \\ \hline 
$0$ & ${{n} \choose 0}$ & ${{n} \choose 1}$  & ${{n} \choose 2} $ & $\cdots$ & ${{n} \choose {n-1}}$ &  $b_{n+1,n+2} = {n \choose n}$ \\ \hline 
\end{tabular} \] 

When $d$ is arbitrary and $\dim_k(R_1) = 1$ the shape is analogous:
\[ \begin{tabular}{|c|c|c|c|c|c|c|c|} 
\hline $b_{0,0} = {{n} \choose 0}$ & ${{n} \choose 1}$ & ${{n} \choose 2}$  & ${{n} \choose 3} $ & $\cdots$ & ${{n} \choose {n}}$ & $0$  \\ \hline 
$0$ & $0$ & $0$ & $0$  & $0$ & $0$ &  $0$  \\ \hline 
 $\vdots$ & $\vdots$ & $\vdots $ & $\vdots$  & $\vdots$ & $\vdots$ &  $\vdots$  \\ \hline 
 $0$ & $0$ & $0$ & $0$  & $0$ & $0$ &  $0$  \\ \hline 
 $0$ & ${{n} \choose 0}$ & ${{n} \choose 1}$  & ${{n} \choose 2} $ & $\cdots$ & ${{n} \choose {n-1}}$ &  $b_{n+1,n+1+d} = {n \choose n}$ \\ \hline 
\end{tabular} \] 

\medskip

\nt as the resolution of the ideal $I_1$ and dual account for the betti table. 

\medskip

Note that whenever $\dim_k(R_1) < n+1$, the Gorenstein ring $R_\bullet$ is the quotient:
$\Sym^\bullet(S_1/I_1) \rightarrow R_\bullet$
of a symmetric algebra in fewer variables, so the only ``new'' quotients are those for which 
$I_1 = 0$ and $I_\bullet$ is generated in degree two or more. 

\medskip

\nt {\bf Example 2.6.} When $d = 2$ and $\dim_k(R_1) = n+1$, then: 
\[ \dim_k(R_2) = 1 \ \mbox{and} \ \dim_k(I_2) = {{n+2} \choose 2} - 1 = \frac{n^2 + 3n}{2} \] 
 which completely determine the rest of the betti table. Thus:

\[ \begin{tabular}{|c|c|c|} \hline
 1 & 0 & 0 \\ \hline
 0 & 2 & 0 \\ \hline
 0 & 0 & 1 \\ \hline  \end{tabular} \]

 \[ \begin{tabular}{|c|c|c|c|} \hline
 1 & 0 & 0 & 0  \\ \hline
 0 & 5 & 5 & 0 \\ \hline
 0 & 0 & 0 & 1 \\ \hline  \end{tabular} \] 

 \[ \begin{tabular}{|c|c|c|c|c|} \hline
 1 & 0 & 0 & 0 & 0  \\ \hline
 0 & 9 & 16 & 9 & 0 \\ \hline
 0 & 0 & 0 &  0 & 1 \\ \hline  \end{tabular} \] 

\nt are the tables for non-degenerate quadratic forms when $n = 1,2,3$. 

\medskip

\nt {\bf Remark.} In the second table, duality implies that the map: 
 \[ f: S(-2)^5 \rightarrow S(-3)^5 \ \mbox{is a skew symmetric matrix of linear forms}  \]  

\section{Exact Complexes on $\PP^n, d= 3$ and $n = 1$} 

The sheafification of the  free-module resolution of $I_\bullet$ is a complex:
\[ (*) \ 0 \rightarrow E_{n+1} \rightarrow \cdots \rightarrow E_1 \rightarrow E_0 = {\mathcal O}_{\PP^n} \rightarrow 0 \] 
of direct sums of line bundles, where:
\[ E_i = \bigoplus {\mathcal O}_{\PP^n}(-j)^{b_{ij}} \] 
and in particular, by duality, 
$E_{n+1} = \omega_{\PP^n}(-d)$.

\medskip 

This provides constraints on the entries of the betti table via the Chern characters of the entries of the complex 
of vector bundles since:
\[ \sum _{i = 0}^{n+1} (-1)^i \ch(E_i) = 0 \] 
which is actually $n+1$ equations (one in each degree):

\medskip

\nt {\bf Constraints on the $b_{i,j}$:} For each $e = 0,...,n$ 
\[ \sum_{i,j} (-1)^{i} b_{i,j} j^e = 0 \]     

When $d = 3$, however, this does not determine the betti table even when $I_1 = 0$. 

\medskip

\nt {\bf Example 3.1.} When $n = 1, d = 3$ (and $\dim(R_\bullet) = (1, 2, 2, 1)$) the betti table is:

\[ \begin{tabular}{|c|c|c|c|c|} \hline
 1 & 0 & 0   \\ \hline
 0 & 1 & 0  \\ \hline
 0 & 1 & 0   \\ \hline  
 0 & 0 & 1  \\ \hline
 \end{tabular} \] 
as the Gorenstein rings in this case have complete intersection resolutions:
\[ 0 \rightarrow {\mathcal O}_{\PP^1}(-5) \rightarrow {\mathcal O}_{\PP^1}(-2) \oplus  {\mathcal O}_{\PP^1}(-3) \rightarrow {\mathcal O}_{\PP^1}  \rightarrow 0 \]

\nt {\bf Example 3.2.} When $n = 2, d = 3$ (and $\dim(R_\bullet) = (1, 3, 3, 1)$) the betti table is:

\[ \begin{tabular}{|c|c|c|c|c|} \hline
 1 & 0 & 0 & 0   \\ \hline
 0 & 3 & $b$ & 0  \\ \hline
 0 & $b$ & 3 &  0  \\ \hline  
 0 & 0 & 0 & 1 \\ \hline
 \end{tabular} \] 
and we only conclude from the constraints that $b_{1,3} = b = b_{2,3}$. Actually, we can say more when we analyze things a bit. Since the map $f$ of vector bundles in the complex below is 
skew symmetric and has co-rank one: 
\[ \begin{array}{ccccccccccc} 
& & & &  {\mathcal O}_{\PP^2}(-3)^b & & {\mathcal O}_{\PP^2}(-2)^3  \\
0 & \rightarrow & {\mathcal O}_{\PP^2}(-6) & \rightarrow  & \oplus &   \stackrel f  \rightarrow & \oplus
 &  \rightarrow &  {\mathcal O}_{\PP^2}&  \rightarrow & 0 \\
& & & & {\mathcal O}_{\PP^2}(-4)^3 & &  {\mathcal O}_{\PP^2}(-3)^b  \end{array} \] 
it follows that $b$ is even, and indeed both $b = 0$ and $b = 2$ do occur, the latter for socles $\sigma$ that lie in the (closure of) the variety of 
planes spanned by length three non-collinear subschemes $Z \subset \PP^2$. 

\medskip

\nt {\bf Example 3.3.} When $n = 3$ and $d = 3$ (and $\dim(R_\bullet) = (1, 4, 4, 1)$), the table is
\[ \begin{tabular}{|c|c|c|c|c|c|} \hline
 1 & 0 & 0 & 0  & 0 \\ \hline
 0 & 6 & $b + 5$ & $b$ & 0  \\ \hline
 0 & $b$ & $b+5$  &  6 & 0  \\ \hline  
 0 & 0 & 0 & 0 &  1 \\ \hline
 \end{tabular} \] 

\medskip

\nt and it is a challenge to determine the  possible values of $b$ (and the corresponding loci of socles). Note that as in the previous example, we may consider the span of 
four non-coplanar points whose ideal resolves as:
\[ 0 \rightarrow {\mathcal O}_{\PP^3}(-4)^3 \rightarrow {\mathcal O}_{\PP^3}(-3)^8 \rightarrow 
{\mathcal O}_{\PP^3}(-2)^6 \rightarrow {\mathcal I}_Z \rightarrow 0 \] 
to deduce that $b = 3$ is obtained by socles contained in the span of four points. 
  
\medskip

\nt {\bf Remark.} These examples are related to the betti tables for canonically embedded 
non-hyperelliptic curves of genus $4,5$ and $6$ since the generic restriction of the minimal free resolution of 
a canonically embedded curve to a codimension two subspace is the resolution of a Gorenstein ring with socle in degree $3$. See e.g. \cite{Eis}

\medskip

We may reinterpret the projective space $\PP^n_d$ via Serre duality as: 
\[ \PP^n_d =  \PP(\mbox{H}^0(\PP^n, {\mathcal O}_{\PP^n}(d))) = \PP(\Ext^n({\mathcal O}_{\PP^n}, \omega_{\PP^n}(-d)^\vee)) \] 
parametrizing non-zero morphisms
\[ \sigma \in \mbox{Hom}({\mathcal O}_{\PP^n}, \omega_{\PP^n}(-d)[n] ) \] 
in ${\mathcal D}^b(\PP^n)$, and then we may reinterpret the resolution as a distinguished triangle:
\[ \rightarrow [E_n \rightarrow \cdots \rightarrow E_1] \rightarrow {\mathcal O}_{\PP^n} \stackrel \sigma \rightarrow 
\omega_{\PP^n}(-d)[n] \rightarrow C_\sigma \rightarrow \cdots \]
or, more suggestively, in the form: 
\[ \rightarrow \omega_{\PP^n}(-d)[n-1] \rightarrow  [E_n \rightarrow \cdots \rightarrow E_1]  \rightarrow {\mathcal O}_{\PP^n} \rightarrow  \] 
which will then be a {\bf short exact sequence} in any abelian subcategory of ${\mathcal D}^b(\PP^n)$ that contains both ${\mathcal O}_{\PP^n}$ and 
$\omega_{\PP^n}(-d)[n-1]$ and for which:
\[ \Ext^1_{\mathcal A}({\mathcal O}_{\PP^n}, \omega_{\PP^n}(-d)[n-1]) = \Ext^n_{Coh(\PP^n)}({\mathcal O}_{\PP^n}, \omega_{\PP^n}(-d))  \]

\medskip

The way to proceed when $n = 1$ was explained to each of the authors as graduate students by 
their mutual advisor Robert Lazarsfeld. 
In that case, the embedding: 
\[ \PP^1 \hookrightarrow \PP^1_d = \PP^{d} \] 
is the {\bf rational normal curve}, and the suggestion is realized with ${\mathcal A} = {\mathcal Coh}(\PP^1)$. 
Since every vector bundle on $\PP^1$ is a sum of line bundles, each short exact sequence associated to a socle has the form: 
\[ \sigma: \ 0 \rightarrow {\mathcal O}_{\PP^1}(-2-d) \rightarrow {\mathcal O}_{\PP^1}(-a) \oplus {\mathcal O}_{\PP^1}(-b) \rightarrow {\mathcal O}_{\PP^1} \rightarrow 0 \] 
with
$a + b = d + 2$ (and we assume $a \le b$). From this point of view, the strata of $\PP^{d}$ are jumping loci for
$E(\sigma) = {\mathcal O}_{\PP^1} (-a) \oplus {\mathcal O}_{\PP^1} (-b)$ 
and then if $a < b$, the following diagram is uniquely determined by $\sigma$ 
\[ \begin{array}{cccccccccc} 
& & & & {\mathcal O}_{\PP^1}(-a) & = & {\mathcal O}_{\PP^1}(-a) \\
& & & & \downarrow& & \downarrow  \\
\sigma: \ 0 & \rightarrow & \omega_{\PP^1}(-d)  & \rightarrow & {\mathcal O}_{\PP^1}(-a) \oplus {\mathcal O}_{\PP^1}(-b) & \rightarrow & {\mathcal O}_{\PP^1} & \rightarrow & 0 \\
& & || & & \downarrow & & \downarrow \\ 
s: \ 0 & \rightarrow &  \omega_{\PP^1}(-d) & \rightarrow  & {\mathcal O}_{\PP^1}(-b) & \rightarrow & {\mathcal O}_Z & \rightarrow & 0 \\
\end{array} \] 

Each integer $0 < a < \frac{d+2}2$ indexes a stratum in $\PP^d$ by:

\medskip

$\bullet$ a choice of quotient sheaf ${\mathcal O}_Z$ for $l(Z) = a$, i.e. a point of $\Sym^a(\PP^1)$

\medskip

$\bullet$ a choice of extension (modulo scalars): 
\[ s \in \mbox{Ext}^1({\mathcal O}_Z,\omega_{\PP^1}(-d)) \] 
that is generic, in the sense that the associated object is the line bundle ${\mathcal O}_{\PP^1}(-b)$. 
This space is the linear span of the subscheme $Z\subset \PP^1 \subset \PP^d$ via:
\[ \Ext^1_{\PP^1}({\mathcal O}_Z,\omega_{\PP^1}(-d)) = \mbox{Hom}({\mathcal O}_{\PP^1}(-d), {\mathcal O}_Z)^\vee = \mbox{H}^0(\PP^1, {\mathcal O}_Z(d))^\vee  \]  
and the generic extensions are the complement of the spans of subschemes $Z' \subset Z$. 

\medskip

Thus we obtain open subsets
$U_a \subset \PP_a \ \mbox{of projective bundles of abstract spans $|Z|$}$ 
embedding into $\PP^d$ as locally closed subschemes making up a {\it stratification} with
\[ Z_a = \overline {U_a} = U_a \cup U_{a-1} \cup \cdots \cup U_1  \] 

When $d$ is odd, this is the full story. When $d$ is even, the generic stratum indexes: 
\[ \sigma: \ 0 \rightarrow {\mathcal O}_{\PP^1}(-2b) \rightarrow {\mathcal O}_{\PP^1}(-b)^2 \rightarrow {\mathcal O}_{\PP^1} \rightarrow 0 \] 

This picture generalizes when $\PP^1$ is replaced by a curve $C$ of any genus. Vector bundles on curves
do not split in general, but the unstable ones do filter. Let:
\[ \PP^{d-g} = \PP(\mbox{H}^0(C, L)) = \PP(\Ext^1({\mathcal O}_C, \omega \otimes L^{-1})^\vee) \] 
for a fixed line bundle $L$ on $C$ of degree $d \ge 2g - 1$. Then a ``generalized socle''
\[ \sigma: 0 \rightarrow \omega_C \otimes L^{-1} \rightarrow E \rightarrow {\mathcal O}_C \rightarrow 0 \] 
determines the rank two vector bundle $E$ which is either:

\medskip

(a) unstable, in which case $A^\vee \subset E$ is a uniquely determined sub-line bundle (with quotient line bundle!) 
of maximal degree $-a$ where $a < ((2 - 2g) + d)/2$ 

\medskip

(b) semi-stable, in which case $\sigma$ is in the regular locus of the rational map
\[ f: \PP^{d-g} --> {\mathcal M}_C(2,L) \] 
to the projective moduli space of semi-stable bundles. 

\medskip

In case (a), we have the same diagram as for $\PP^1$: 

\[ \begin{array}{cccccccccc} 
& & & & A^\vee & = & A^\vee \\
& & & & \downarrow& & \downarrow  \\
\sigma: \ 0 & \rightarrow & \omega_C \otimes L^{-1}  & \rightarrow & E & \rightarrow & {\mathcal O}_C & \rightarrow & 0 \\
& & || & & \downarrow & & \downarrow \\ 
s: \ 0 & \rightarrow &   \omega_C \otimes L^{-1}  & \rightarrow  &  B^\vee & \rightarrow & {\mathcal O}_Z & \rightarrow & 0 \\
\end{array} \] 

\medskip

\nt and the indeterminacy of the map to moduli is stratified by locally closed sets: 
\[ U_a \subset \PP^{d - g}  \] 
isomorphic to open subsets of the projective bundles over Sym$^a(C)$  with fibers: 
\[ \PP(\mbox{H}^0(C, L_Z))  \] 
of spans of the schemes $Z\subset C \subset \PP^{d-g}$ of length $a$.

\medskip

This has proven to have many useful applications, and it is this construction that we wish to reinterpret and generalize to  higher-dimensional projective spaces, 
starting with the projective plane $\PP^2$.

\section{Stabilities on $\PP^n$} 

Let ${\mathcal D}$ be the bounded derived category of coherent sheaves on $\PP^n$.

\medskip

\nt {\bf Definition 4.1.} A  {\bf stability condition} on ${\mathcal D}$ is:

\medskip

$\bullet$ a bounded $t$-structure with (abelian) heart ${\mathcal A}$

\medskip

$\bullet$ a linear ``central charge'' defined on the $K$-group of vector bundles on $\PP^n$: 
\[ Z: K(\PP^n)  \rightarrow \CC \ \mbox{such that}  \] 

(i) $0 < \arg(Z(A)) \le 1$ for all objects $A$ of ${\mathcal A}$, where $z = re^{\pi i \arg(z)}$, and we assign:
\[ \arg(A[m]) := \arg(A) + m \ \mbox{for shifts by $m \in \ZZ$ of objects of ${\mathcal A}$} \]

(ii) the stable objects of ${\mathcal D}$ are shifts $E[m]$ of objects $E$ of ${\mathcal A}$ satisfying
\[ \arg(Z(A)) < \arg(Z(E)) \ \mbox{whenever} \   A \subset E \] 

(iii) the semi-stable objects of ${\mathcal D}$ are shifts of objects $F$ of  ${\mathcal A}$ that filter:  
\[ 0 = F_0 \subset F_1 \subset \cdots \subset F_{n(F)} = F \ \mbox{with} \ E_i = F_i/F_{i-1} \ \mbox{stable with fixed $\arg$} \]

\medskip

(iv) every object $C$ of ${\mathcal D}$ filters in the derived sense via a finite sequence: 
\[ C_{i-1} \stackrel {f_i} \rightarrow C_{i} ; \ C_0 = 0 \ \mbox{and} \ C_{m(C)} = C \] 
of morphisms with the property that the {\bf cone} $C_{f_i}$ of each morphism satisfies:
\[ C_{f_i} \ \mbox{is semistable and the} \ r_i := \arg(C_{f_i} ) \in \RR, \ \mbox{satisfy} \ r_1 > r_2 > ... > r_{m(A)} \] 

If $C = A$ is an object of ${\mathcal A}$, then the filtration consists of injective morphisms: 
\[ 0 = A_0 \subset A_1 \subset \cdots \subset A_{m(A)} = A \] 
with $F_i = A_i/A_{i-1}$ and $1 \ge \arg(F_1) > \cdots > \arg(F_m(A)) > 0$.

\medskip

Such a sequence is (essentially) unique and is called the {\bf Harder-Narasimhan} filtration of the 
object $C$ with respect to the stability condition $({\mathcal A}, Z)$. 

\medskip

\nt {\bf Note.} The {\it finiteness} (of each $n(F)$ and $m(A)$)  is the issue, 
once a bounded $t$-structure ``supporting'' the stability condition has been specified. Bridgeland showed \cite{Bay} that 
stability conditions deform along with deformations of the linear map $Z$ when it has the {\it support} property: 

\medskip

(v) There is a quadratic form $Q$ on $\RR^{n+1} = K({\mathcal D}) \otimes \RR$ such that:

\medskip

$\bullet$ $Q(Z(F)) \ge 0$ for all semi-stable objects, and 

\medskip

$\bullet$ $Q < 0$ on the kernel of $Z_\RR: \RR^{n+1} \rightarrow \CC = \RR^2$. 

\medskip

One property of stability conditions that we will need is:

\medskip

\nt {\bf Rotation.} An action of $\RR$ on the locus Stab$({\mathcal D})$ of stability conditions is given by: 
\[ r\cdot ({\mathcal A}, Z) = ({\mathcal A}^\#, Z^\# = e^{-\pi i r}\cdot Z) \ \mbox{for} \ r = n + \theta \ \mbox{with} \ n = \lfloor r \rfloor \] 
where each object of $A$ of ${\mathcal A}^\#$ has a Harder-Narasimhan filtration  
divided into: 
\[ 0 \subset A_0 \subset \cdots \subset A_{m_1(A)} \ \mbox{and} \ A_{m_1(A)} \subset \cdots \subset A_{m_1(A) + m_2(A)} = A  \] 
where the first part of the filtration satisfies $A_i/A_{i-1} = T_i[n]$ for semistable $T_i \in {\mathcal A}$ and strictly decreasing ``large arguments'' $\theta_i := \arg(T_i) > \theta$, and
the second part 
satisfies with $A_j/A_{j-1} = F_j$ and strictly decreasing ``small args'' $\theta \ge \theta_j = \arg(F_j)$
(augmented by the extra shift to become the large args in ${\mathcal A}^\#$). 

\medskip

\nt {\bf Example 4.1.} Consider the Mumford stability condition on $\PP^1$ with: 

\medskip

$\bullet$ ${\mathcal A} = Coh(\PP^1)$ the heart of the standard $t$-structure on ${\mathcal D}^b(\PP^1)$ and 

\medskip

$\bullet$ $Z({\mathcal O}_{\PP^1}) = (0,1)$ and $Z(\CC_p) = (-1,0)$ 

\medskip

\nt with
$Z(E) = (-\ch_1(E), \ch_0(E))$ lying on the ray $\arg(z) = 1$ when $E$ is a torsion sheaf and the strict upper 
half plane otherwise. The only stable objects are the skyscraper sheaves $\CC_p$ and line bundles ${\mathcal O}_{\PP^1}(d)$ (and all shifts). The 
semi-stables are torsion sheaves and direct sums $V \otimes_k {\mathcal O}_{\PP^1}(d)$ (and their shifts). 

\medskip

When rotated by $r = \frac 12$, the new heart ${\mathcal A}^\#$ contains:
\[ {\mathcal O}_{\PP^1}[1] \ \mbox{and} \ {\mathcal O}_{\PP^1}(1) \] 
with $\arg(Z^\#({\mathcal O}_{\PP^1}[1])) = 1$ and  $\arg(Z^\#({\mathcal O}_{\PP^1}(1))) = \frac 14$ and  
it follows that ${\mathcal A}^\#$ is the category of morphisms (representations of the Kronecker quiver): 
\[ V_1 \otimes {\mathcal O}_{\PP^1} \rightarrow V_0\otimes {\mathcal O}_{\PP^1}(1) \] 
Note that the semi-stable objects are the same 
as for the stability condition on ${\mathcal A}$.  

\medskip

Under the automorphism $E \mapsto E\otimes {\mathcal O}_{\PP^1}(d)$ of the derived category ${\mathcal D}^b(\PP^1)$  the subcategory $Coh(\PP^1)$ remains fixed 
 but the 
categories ${\mathcal A}^\#$ are translated, replacing ${\mathcal O}_{\PP^1}[1]$ and ${\mathcal O_{\PP^1}}(1)$ with 
${\mathcal O}_{\PP^1}(-d)[1]$ and ${\mathcal O}_{\PP^1}(-d+1)$, respectively. In particular, if $d = 1$, then 
the Mumford stability condition on ${\mathcal A} = Coh(\PP^1)$ is replaced with an equivalent stability condition setting: 
\[ Z({\mathcal O}_{\PP^1}) = (1,1) \ \mbox{and} \ Z(\CC_p) = (-1,0) \] 

This is the precise stability condition that we generalize from $\PP^1$ to  $\PP^n$.

\medskip

The stability condition we will use for $n > 1$ and $d = 2e$ is obtained in two steps.

\medskip

\nt {\bf Step 1.} Let ${\mathcal C}$ be the abelian category of complexes of vector bundles  of the form: 
\[ {\mathcal C} = \{ {V_n \otimes \mathcal O}_{\PP^n}(-n) \rightarrow \cdots \rightarrow V_0 \otimes {\mathcal O}_{\PP^n} \} \] 
This is the heart of one of many ``quivery''  Beilinson $t$-structures on ${\mathcal D}$ that supports a {\it complex manifold} of 
stability conditions with: 
\[ Z({\mathcal O}_{\PP^n}(-i)[i]) =: \zeta_i \]  
for {\bf any} $\zeta_0,...,\zeta_n$ satisfying $0 < \arg(\zeta_i) \le 1$. The moduli of stable objects
for such stability conditions is obtained through GIT, as first shown by Alastair King \cite{Ki}.

\medskip

Building on the $\QQ[t]$-valued linear {\bf Hilbert polynomial}
\[ \chi(E(t)) := \sum_{j = 0}^n (-1)^j \dim(\mbox H^j(\PP^n, E\otimes {\mathcal O}_{\PP^n}(t))) \] 
we choose a stability condition by setting: 
\[ Z(E) = (\chi'(E),   \chi(E)) := (\chi'(E(t))|_{t = 0} , \chi(E(t))|_{t = 0} ) \] 
noting that: 
\[ Z({\mathcal O}_{\PP^n}) = \left(\sum_{i=1}^n \frac 1i , 1\right) \ \mbox{and} \  \] 
\[ \arg(Z({\mathcal O}_{\PP^n}(-i)[i])) = 1 \ \mbox{for all} \ i = 1,...,n \]  

\medskip

\nt so  that the objects of ${\mathcal C}$ all satisfy $0 < \arg(Z(A)) \le 1$ as required. 

\medskip

This is ${\mathcal A}^\#$ for the modified Mumford slope when $n = 1$, and  then in general,

\medskip

\nt {\bf Proposition 4.1.} Skyscraper sheaves $\CC_p$  and all line bundles ${\mathcal O}_{\PP^n}(d)$ are stable. 

\medskip

{\bf Proof.} Stability of the line bundles ${\mathcal O}_{\PP^n}(-i)$ for $i = 0,...,n$ follows from the fact that their shifts
${\mathcal O}_{\PP^n}(-i)[i]$ are {\bf simple} objects of ${\mathcal C}$ (no non-zero sub-complexes).  
The skyscraper sheaves are objects of ${\mathcal C}$ via the Koszul complex:
\[ \CC_p = [  \wedge^nW \otimes {\mathcal O}_{\PP^n}(-n) 
\rightarrow \cdots \rightarrow W\otimes {\mathcal O}_{\PP^n}(-1) \rightarrow  k\otimes {\mathcal O}_{\PP^n}] \]
with $W = \mbox{H}^0(\PP^n, {\mathcal I}_p(1))$ as discussed earlier. They are stable because $V_0 = k$ 
and as a result
any sub-complex $E$ with dimension vector $(d_i = \dim(V_i))$ either satisfies: 

\medskip

$\bullet$ $d_0 = 1$ and all $d_i \le {n \choose i}$ and some $d_i < {n \choose i}$, in which case:  
\[ Re(Z(E))  >  0 \ \mbox{and} \ \arg(Z(E)) < \arg(Z(\CC_p)) = \frac 12   \] 
or else 

\medskip

$\bullet$ $d_0 = 0$. But for $i > 0$,
\[ \mbox{Hom}({\mathcal O}_{\PP^n}(-i)[i], \CC_p) = \Ext^{-i} ({\mathcal O}_{\PP^n}(-i), \CC_p) = 0 \]
and so by induction, Hom$(E,\CC_p) = 0$  if $V_0 = 0$. So there are no such sub-complexes.  

\medskip

For the line bundles ${\mathcal O}_{\PP^n}(d)$ and $d > 0$, we note that\footnote{This was pointed out to us by Marin Petkovic.}:

\medskip

$\bullet$ both $Z$ and line bundles ${\mathcal O}_{\PP^n}(d)$ are invariant under the action of GL$(n,\CC)$ 

\medskip

$\bullet$ the Eagon-Northcott complex of GL$(n,\CC)$-representations: 
\[ {\mathcal O}_{\PP^n}(d) = [ EN_n\otimes {\mathcal O}_{\PP^n}(-n) \rightarrow \cdots \rightarrow EN_0 \otimes  {\mathcal O}_{\PP^n}] \ \mbox{with} \ 
R_0 = \Sym^dV  \]  
commences with an {\it irreducible} representation of GL$(V)$. For this reason, every  {\bf maximal destabilizing complex} 
(first term in the Harder-Narasimhan filtration) is a complex of sub-representations, and then by the irreducibility of 
the representation Sym$^d(\CC^n)$, the argument for the skyscraper sheaves also applies here. 

\medskip

The remaining shifts $\omega_{\PP^n}(-d)[n]$ for $d \ge 0$ lying in ${\mathcal C}$ are also stable thanks to the following 
crucial feature of the central charge $Z(E) = (\chi'(E), \chi(E))$: 

\medskip

The category ${\mathcal C}$ is invariant under the involution of ${\mathcal D}$:
\[ E \mapsto E^\vee  := RHom_{\mathcal D} (E, {\mathcal O}_{\PP^n}(-n)) [n] \] 
and the central charge satisfies:
\[ Z(E^\vee(-1)[n]) = (-\chi'(E), \chi(E)) \] 

This follows from Serre duality, as $E^\vee (-1) = RHom(E,\omega_{\PP^n})$ and
\[ \chi(RHom(E,\omega_{\PP^n})[n] (t)) =  \chi(E(-t))  \] 

Thus, not only is it the case that the $\omega_{\PP^n}(-d)[n]$ are stable, but:
\[ Z(\omega_{\PP^n}(-d)[n])  \ \mbox{is the reflection of} \ Z({\mathcal O}_{\PP^n}(d))  \ \mbox{across the imaginary axis}\] 

This leads us to the definition of ${\mathcal A}$ generalizing the $n = 1$ case.

\medskip

{\bf Step 2.} ${\mathcal A} = {\mathcal C}^\#$ is the rotation of the stability condition $({\mathcal C}, Z)$ by $r = -\frac 12$. 

\section{The Even Degree Case $d = 2e$} 

Recall from Example 2.2 that the socles in degree two  are symmetric maps: 
\[ q:V \rightarrow V^\vee \ \mbox{where} \ V = (S_1)^\vee \] 
and that the projective space of symmetric matrices stratifies by {\bf rank}: 
\[ \PP^n_2 = \coprod U_r  \] 
where
\[ \begin{array}{cccccc}  \PP(\Sym^2(Q))^0  & \rightarrow &U_r &  \subset &  \PP(\Sym^2({\mathcal Q}))  \\
\downarrow & & \downarrow  & \swarrow \\ Q & \in & \mbox{Gr}(V,r) \end{array} \] 
i.e. $U_r$ is the fiber bundle over the Grassmannian of $r$-dimensional quotients of $V$ consisting of injective symmetric maps of the quotient space, generating:
\[ (\dagger) \ V \stackrel \lambda \rightarrow Q \rightarrow Q^\vee \stackrel {\lambda^\vee}\rightarrow V^\vee \] 

\medskip

We want to generalize this to higher even values of $d$ by reinterpreting the socle as a 
morphism of ``balanced'' stable objects in ${\mathcal C}$:
\[ \sigma: {\mathcal O}_{\PP^n}(e) \rightarrow \omega_{\PP^n}(-e)[n] \] 
and then mimic the strategy above by finding: 

\medskip

(a) the full set of ``ranks'' (multi-degrees) $r$ of maps $\sigma$, and for each $r$, finding

\medskip

(b) the {\bf Hilbert scheme} of quotient complexes 
$\Hilb_{{\mathcal C}} ({\mathcal O}_{\PP^n}(e), r)$, and 

\medskip

(c) the open subset $U_r$ of the family:
\[ \PP(\Hom_{\mathcal C}({\mathcal O}_\Lambda(e) ,\omega_{\PP^n}(-e)[n]) \]  
so that the composition ${\mathcal O}_{\PP^n}(e) \rightarrow \omega_{\PP^n}(-e)[n]$ has full rank ${r-1}$.

\medskip

In this way, the locally closed loci subsets $U_r$ of symmetric matrices are:

\[ \begin{array}{cccccc}  \PP(\Hom({\mathcal O}_\Lambda(1) ,\omega_{\PP^n}(-1)[n]))^0  & \rightarrow &U_r &  \subset &  \PP(\Hom({\mathcal Q},\omega_{\PP^n}(-1)[n]))  \\
\downarrow & & \downarrow  & \swarrow \\ Q & \in & \mbox{Hilb}({\mathcal O}_{\PP^n}(1),r)  \end{array} \] 
and the Grassmannian is reinterpreted as the Hilbert scheme of quotients: 
\[ \ {\mathcal O}_{\PP^n}(1) \rightarrow {\mathcal O}_\Lambda(1) \ \mbox{for linear subspaces} \ \PP^{r-1} \cong \Lambda \subset \PP^{n}  \] 

 In this case, 
the factorization of a map to $\omega_{\PP^n}(-1)[n]$ through:
\[ {\mathcal O}_{\Lambda}(1) \rightarrow \omega_\Lambda(-1)[r] \rightarrow \omega_{\PP^n}(-1)[n] \] 
generating the socle/symmetric matrix ``on $\Lambda$'' is automatic, following from: 
\[ \Ext^n({\mathcal O}_\Lambda, \omega_{\PP^n}(-1)) = \Hom({\mathcal O}_{\PP^n}(-1), {\mathcal O}_{\Lambda}(1))^\vee  = 
\Hom({\mathcal O}_{\Lambda}(-1), {\mathcal O}_{\Lambda}(1))^\vee \] 
and another application of Serre duality on $\Lambda$, so that: 
\[  {\mathcal O}_{\PP^n}(1) \stackrel \lambda \rightarrow  {\mathcal O}_\Lambda(1) \rightarrow  \omega_{\PP^n}(-1)[n]  \] 
is the socle version of the decomposition $(\dagger)$ above.

\medskip

The dimension one case fits this reinterpretation for all $e$, via the complexes 

\[ \sigma: [{\mathcal O}_{\PP^1}(-1)^e \rightarrow {\mathcal O}_{\PP^1}^{e+1}]  \rightarrow  [{\mathcal O}_{\PP^1}(-1)^{e+2} \rightarrow {\mathcal O}_{\PP^1}^{e+1}]  
 \] 

\nt and

\medskip

($0$) the ``full rank'' maps are injective with a common cokernel ${\mathcal O}_{\PP^1}(-1)^2[1]$. 

\medskip

($r$) otherwise, the socle factors uniquely through a surjection: 
\[ {\mathcal O}_{\PP^1}(e) \rightarrow {\mathcal O}_Z(e) = [{\mathcal O}_{\PP^1}^r \rightarrow {\mathcal O}_{\PP^1}^r] \]  
for a unique closed subscheme $Z \subset \PP^1$ of length $r$, 
followed by an injection:
\[ \overline \sigma: {\mathcal O}_Z(e) = [{\mathcal O}_{\PP^1}(-1)^r \rightarrow  {\mathcal O}_{\PP^1}^r] \rightarrow 
[{\mathcal O}_{\PP^1}(-1)^{e+2} \rightarrow {\mathcal O}_{\PP^1}^{e+1}]   \] 
yielding an exact sequence:
\[ \ker(\sigma) = {\mathcal O}_{\PP^1}(e - r) \rightarrow {\mathcal O}_{\PP^1}(e) \stackrel \sigma \rightarrow 
\omega_{\PP^1}(-e)[1] \rightarrow  \omega_{\PP^1}(-e+r)[1] = \coker(\sigma) \] 
and the familiar consolidated short exact sequence of vector bundles:  
\[ \sigma: \ 0 \rightarrow \omega_{\PP^1}(-e) \rightarrow {\mathcal O}_{\PP^1} (e - a) \oplus \omega_{\PP^1}(-e+a) \rightarrow {\mathcal O}_{\PP^1}(e) \rightarrow 0 \] 
The splitting of the internal term follows here from the fact that:
\[ \Ext^2(\omega_{\PP^1}(-e+a)[1], {\mathcal O}_{\PP^1}(e + a)) = \Hom({\mathcal O}_{\PP^1}(e+a), {\mathcal O}_{\PP^1}(e-a)) = 0 \] 

\medskip

The interpretation of $U_r$ as a fiber bundle is, as in the earlier construction: 

\[ \begin{array}{cccccc}  \PP(\Ext^1({\mathcal O}_Z(e) ,\omega_{\PP^1}(-e)))^0  & \rightarrow &U_r &  \subset &  \PP(\Hom({\mathcal Q},\omega_{\PP^n}(-e)[1]))  \\
\downarrow & & \downarrow  & \swarrow \\ {\mathcal O}_Z(e)  & \in & \mbox{Hilb}({\mathcal O}_{\PP^1}(e),r)  \end{array} \] 

\medskip

\nt and then there are two things to verify in this Example to be sure the stratification covers $\PP^d$, namely that if $\sigma$ factors through a quotient {\bf complex}
\[ \sigma: {\mathcal O}_{\PP^1}(e) \rightarrow Q \stackrel{ \overline \sigma} \hookrightarrow \omega_{\PP^1}(-e)[1] \] 
then  $Q$ is a quotient {\bf coherent sheaf}, which must therefore be an ${\mathcal O}_Z(e)$ since: 
\[ [{\mathcal O}_{\PP^1}(-1)^e \rightarrow {\mathcal O}_{\PP^1}^{e+1}] \rightarrow [{\mathcal O}_{\PP^1}(-1)^r\rightarrow {\mathcal O}_{\PP^1}^r] \]
can visibly be surjective only when $e \ge r$. 

\medskip

For this last remaining point and to tackle the cases $n \ge 2$, we rely on: 

\medskip

\nt {\bf Theorem 5.0.}  For the stability conditions on ${\mathcal D} = {\mathcal D}^b(\PP^n)$ from \S 4: 
\[ ({\mathcal C},Z) \ \mbox{and} \ ({\mathcal A} = {\mathcal C}^\#, Z^\# = \pi \cdot Z) \]

(a) ${\mathcal O}_{\PP^n}(e)$ and $\omega_{\PP^n}(-e)[n-1]$ are stable objects of ${\mathcal A}$. 

\medskip

(b) The short exact sequence of objects of ${\mathcal A}$ determined by each socle $\sigma$:
\[ \sigma: 0 \rightarrow \omega_{\PP^n}(-e)[n-1]  \rightarrow E(\sigma) \rightarrow {\mathcal O}_{\PP^n}(e) \rightarrow 0 \] 
produces an object $E(\sigma)$ of ${\mathcal A}$ satisfying: 

\medskip

(i) $\arg(Z(E(\sigma))) = 0$ 

\medskip

(ii) $E(\sigma)$ is semi-stable  if and only if $\sigma: {\mathcal O}_{\PP^n}(e) \rightarrow \omega_{\PP^n}(-e)[n]$ is injective in ${\mathcal C}$ 

\medskip

\nt otherwise,

\medskip

(iii) $\sigma$ factors through a coimage $Q(\sigma)$ complex
\[ {\mathcal O}_{\PP^n}(e) \twoheadrightarrow Q(\sigma)  \hookrightarrow \omega_{\PP^n}(-e)[n] \] 
with $H^n(Q(\sigma)) = 0$. 

\medskip

In particular, when $n = 1$, the coimage is a coherent sheaf as noted above.  

\medskip

(iv) The kernel complex $\ker(\sigma)$ satisfies 
\[ H^n(\ker(\sigma)) = H^{n-1}(\ker(\sigma)) = 0 \]
and in particular, when 
$n = 1,2$ it is a coherent sheaf, though in the $n = 2$ case, 
\[ \ker(\sigma) \rightarrow {\mathcal O}_{\PP^2}(e) \]  
may be an inclusion only as complexes, and not as coherent sheaves. 

\medskip

(v) the Harder-Narasimhan filtration of $\ker(\sigma)$ yields a ``maximal destabilizer'' 
\[ F(\sigma) \subset \ker(\sigma) \subset {\mathcal O}_{\PP^n}(e)  \ \mbox{of} \ E(\sigma) \]
as well as the canonically determined rays of the other semi-stable sub-quotients. 

\medskip

{\bf Proof.} In the following diagram,

\medskip

$\bullet$ objects of ${\mathcal C}$ map to the uppper half plane under $Z$, and

\medskip

$\bullet$ objects of ${\mathcal A}$ map to the right half plane under $Z$. 

\bigskip

\begin{center}

\begin{tikzpicture} 

\draw[shift={(0,0)},color=black]  node[below] {$(0,0)$};

\draw[->]  (0,0) -- (2,2);
\draw[shift={(2,2)},color=black] node[right] {${\mathcal O}(e)$};

\draw[->]  (0,0) -- (-2,2);
\draw[shift={(-2,2)},color=black] node[left] {$\omega(-e)[n]$};

\draw[->]  (0,0) -- (2,-2);
\draw[shift={(2,-2)},color=black] node[right] {$\omega(-e)[n-1]$};

\draw[->] (0,0) -- (4,0);
\draw[shift={(4,0)},color=black] node[right] {$E(\sigma)$};

\draw[->] (0,0) -- (.7,.5);
\draw[shift={(.7,.5)},color=black] node[right] {$F(\sigma)$};

\draw[->] (0,0) -- (2,.5);
\draw[shift={(2, .5)},color=black] node[right] {$\ker(\sigma)$};

\end{tikzpicture}  

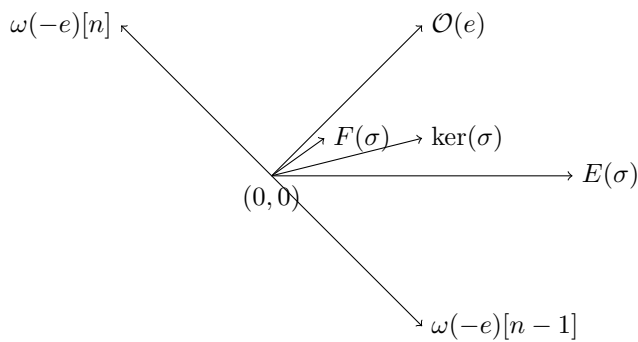
{\mbox Images under the Central Charge Z} 

\end{center}

\bigskip 
(a) By Proposition 1, the objects are $Z$-stable. By direct computation:
\[ \chi({\mathcal O}_{\PP^n}(e)) = { {n + e} \choose n} \ \mbox{and} 
\ \chi'({\mathcal O}_{\PP^n}(e)) = { {n + e} \choose n} \cdot \sum_{i=1}^n \frac 1{e+i} \]
so ${\mathcal O}_{\PP^n}(e)$ maps to the first quadrant (as pictured) and then the objects $\omega_{\PP^n}(-e)[n]$ and 
$\omega_{\PP^n}(-e)[n-1]$ map to quadrants II and IV (also as pictured) by the symmetry. 

\medskip

(bi) follows from the symmetry of Proposition 1 and the additivity of $Z$ with 
\[ Z(E(\sigma)) = (2\chi'({\mathcal O}_{\PP ^n}(e)),0) \] 

(bii) The two ends of the Eagon-Northcott complex:
\[ {\mathcal O}_{\PP^n}(e) = [V_n \otimes {\mathcal O}_{\PP^n}(-n) \rightarrow \cdots \rightarrow V_0 \otimes {\mathcal O}_{\PP^n} ] \] 
have dimensions
\[ \dim(V_0) = { {n + e} \choose n } \ \mbox{and} \ \dim(V_n) = { {n + e - 1} \choose n } \] 
from which it follows that the two ends of the complex for $\omega_{\PP^n}(-e)[n]$ are: 
\[ \dim(W_0) = {  {n + e} \choose n } \ \mbox{and} \ \dim(W_n) = { {n + e +1} \choose n } \] 
and so (as in the case $n = 1$) if $\sigma$ is injective,  the complex for $E(\sigma)[1]$ has the form:
\[ [U_n \otimes {\mathcal O}_{\PP^n}(-n) \rightarrow \cdots \rightarrow U_1 \otimes {\mathcal O}_{\PP^n}(-1)]  \]  
and maximal arg $= 1$ by Proposition 1. It follows from the maximality that $E(\sigma)[1]$, and therefore also $E(\sigma)$, is semi-stable. 
On the other hand, lifting the inclusion: 
\[ \begin{array}{ccccccccc} & & \ker(\sigma) \\
& \swarrow & \downarrow &  \\
E(\sigma) & \rightarrow & {\mathcal O}_{\PP^n}(e) & \stackrel  \sigma \rightarrow & \omega_{\PP^n}(-e)[n] \end{array} \] 
 of the kernel gives a non-zero map to $E(\sigma)$, and since $\ker(\sigma)$ belongs to ${\mathcal C}$, 
 \[ \arg(F_i) > 0 = \arg(E(\sigma)) \] 
for all the Harder-Narsimhan subquotients $F_i$ of $\ker(\sigma)$, and $E(\sigma)$ is unstable.

\medskip

(biii) The cohomology sheaf sequence for the short exact sequence:
\[ 0 \rightarrow Q(\sigma) \rightarrow \omega_{\PP^n}(-e)[n] \rightarrow \coker(\sigma) \rightarrow 0 \] 
commences with:
$0 \rightarrow H^n(Q(\sigma)) \rightarrow \omega_{\PP^n}(-e) \rightarrow H^n(\coker(\sigma)) \rightarrow 
H^{n-1}(Q(\sigma)) \rightarrow 0$ 

\medskip
 
 Keeping in mind that $H^n(C) \subset V_n \otimes {\mathcal O}_{\PP^n}(-n)$ is torsion-free for all 
 complexes, it follows that 
 $H^n(Q(\sigma)) = \omega_{\PP^n}(-e)$ or $0$. But in the former case, the composition:
 \[ \omega_{\PP^n}(-e)[n] = H^n(Q(\sigma))[n]  \rightarrow Q(\sigma) \rightarrow \omega_{\PP^n}(-e)[n] \] 
 is the identity, and $Q(\sigma) \cong \omega_{\PP^n}(-e)[n]$ (via the postulated injective map), which 
 is not possible since the complex ${\mathcal O}_{\PP^n}(e)$ cannot surject onto $\omega_{\PP^n}(-e)[n]$.
  
 \medskip
 
 (biv) From the beginning of the cohomology sequence of complexes when $n \ge 2$: 
 \[ 0 \rightarrow H^n(\ker(\sigma)) \rightarrow 0 \rightarrow H^n(Q(\sigma)) \rightarrow 
 H^{n-1}(\ker(\sigma)) \rightarrow 0 = H^{n-1}({\mathcal O}_{\PP^n}(e)) \ \mbox{of} \ E(\sigma) \] 
 and (biii), we get the desired vanishing. 
 
 \medskip
 
 Note that when $n = 1$ the last term is not $0$ but rather
 ${\mathcal O}_{\PP^n}(e)$ itself. 
 
 \medskip
 
 (bv) follows immediately from the proof of (bii) with $F(\sigma) = F_1$.  
 
 \medskip
 
 \qed
 
We want next to put this into practice in a few examples, especially when $n = 2$. 
We revisit the case  $n = 1$ for one last time for contrast. In this case, 
\[ Z({\mathcal O}_{\PP^1}(-1)[1]) = -1, \ Z(\CC_p) = i \ \mbox{and} \ Z({\mathcal O}_{\PP^1}) = 1 + i \]   
and the diagram of images of $Z$ has a uniform appearance for any value of $e$:

 \begin{center}

\begin{tikzpicture} 

\draw[shift={(0,0)},color=black]  node[below] {$(0,0)$};

\draw[->]  (0,0) -- (-1,0);
\draw[shift={(-1,0)},color=black] node[left] {${\mathcal O}(-1)[1]$};

\draw[->]  (0,0) -- (0,1);
\draw[shift={(0,1)},color=black] node[above] {$\CC_p$};

\draw[->]  (0,0) -- (2,0);
\draw[shift={(2,0)},color=black] node[right] {$E(\sigma)$};

\draw[->]  (0,0) -- (1,1);
\draw[shift={(1,1)},color=black] node[right] {$\mathcal O$};
\draw[shift={(1,1)},color=black] node{$\bullet$};

\draw[->]  (0,0) -- (1,2);
\draw[shift={(1,2)},color=black] node[right] {$\mathcal O(1)$};
\draw[shift={(1,2)},color=black] node{$\bullet$};

\draw[->]  (0,0) -- (1,4);
\draw[shift={(1,4)},color=black] node[right] {$\mathcal O(e)$};

\draw[shift={(1,3)},color=black] node[right] {$\vdots$};
\end{tikzpicture}  

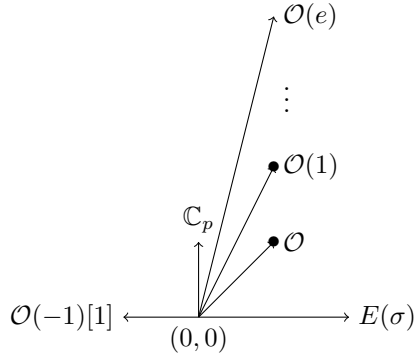
Image of $Z$ for $n = 1$ and $d = 2e$ 

\end{center} 

\medskip

In contrast, each value of $e$ has a different look when $n = 2$. 

\begin{center}

\begin{tikzpicture} 

\draw[shift={(0,0)},color=black]  node[below] {$(0,0)$};

\draw[->]  (0,0) -- (-1,0);
\draw[shift={(-2,0)},color=black] node[above] {${\mathcal O}(-1)[1]$};
\draw[shift={(-2,0)},color=black] node[below] {${\mathcal O}(-2)[2]$};

\draw[->]  (0,0) -- (0,1);
\draw[shift={(0,1)},color=black] node[above] {$\CC_p$};

\draw[->]  (0,0) -- (3,1);
\draw[shift={(3,1)},color=black] node[right] {${\mathcal O}$};
\draw[shift={(3,1)},color=black] node{$\bullet$};

\draw[->]  (0,0) -- (5,3);
\draw[shift={(5,3)},color=black] node[right] {$\mathcal O(1)$};

\draw[->]  (0,0) -- (5,2);
\draw[shift={(5,2)},color=black] node[right] {$\mathcal I_p(1)$};
\draw[shift={(5,2)},color=black] node{$\bullet$};

\draw[->]  (0,0) -- (5,1);
\draw[shift={(5,1)},color=black] node[right]{$\mathcal I_{pq}(1)$};
\draw[shift={(5,1)},color=red] node{$\bullet$};

\draw[shift={(4,2)},color=red] node{$\bullet$};
\draw[shift={(2,1)},color=red] node{$\bullet$};

\end{tikzpicture}  

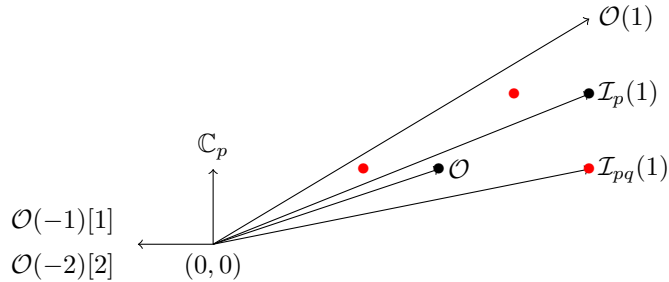
Image of $Z$ for $n = 2$ and $e = 1$ 
\end{center}

The black dots are the $Z$ values of sheaves that {\bf do} arise as kernels and are stable, 
agreeing with the maximal destabilizers. The stratification of:
$\PP^2_2 = \PP^5$
that results matches the stratification by the rank of the quadratic form, as we have seen, with: 
\[ {\mathcal O}_{\PP^2}(1) /{\mathcal O}_{\PP^2} = {\mathcal O}_l(1) \ \mbox{and} 
\ {\mathcal O}_{\PP^2}(1) /I_p(1)  = {\mathcal O}_p(1) \] 
as the two coimages. Then the connection with Gorenstein ring invariants is: 

\medskip

$\bullet$ Open stratum of semi-stable = full rank maps are palindromic:
\[  [1 \ 3 \ 1 ] \]
Hilbert functions, and:
\[ E(\sigma) = [ {\mathcal O}_{\PP^2}(-2)^5 \rightarrow {\mathcal O}_{\PP^1}(-1)^5 ] \] 
is a (skew) map of constant corank one. 

\medskip

$\bullet$ The stratum destabilized by $\ker(\sigma) = {\mathcal O}_{\PP^2}$ has Hilbert function
\[ [1 \ 2 \ 1] \] 
and Harder-Narasimhan filtration of $E(\sigma)$ given by: 
\[ {\mathcal O} = E_1  \subset E_2 \subset E_3 = E(\sigma) \ \mbox{with} \ \]
\[ F_1 = {\mathcal O}_{\PP^2}, F_2 = [{\mathcal O}_{\PP^2}(-2)^2 \rightarrow {\mathcal O}_{\PP^2}(-1)^2], F_3 =  F_1^\vee = {\mathcal O}_{\PP^2}(-3)[1] \] 

\medskip

$\bullet$ The deepest stratum destabilized by $\ker(\sigma) = {\mathcal I}_p(1)$ has Hilbert function
\[ [ 1 \ 1 \ 1] \]  
(these are the points of the Veronese embedding) and filtration:
\[ {\mathcal I}_p(1) = E_1 \subset E_2 = E(\sigma) \ \mbox{with} \] 
\[ F_1 = I_p(1) = [ {\mathcal O}_{\PP^2}(-2) \rightarrow {\mathcal O}_{\PP^2}(-1)^2], F_2 = I_p(1)^\vee \] 

The red dots do not contribute, as they are not kernels even though their $Z$-value 
makes being a kernel a theoretical possibility. This is for two reasons. First, 
$Z(I_{pq}(1))$ lies under the slope of $Z({\mathcal O}_{\PP^2})$ and therefore is not in the 
category ${\mathcal C}$. 

\medskip

The other two red dots do not contribute because the only semi-stable coherent sheaves ${\mathcal F}$ in the category ${\mathcal C}$ with these $Z$-values are {\bf torsion} (rank zero) sheaves, and therefore 
there is no non-zero map from one of these  to ${\mathcal O}_{\PP^2}(1)$. 

\medskip

The red dots proliferate when $e > 1$, so we need to find a systematic method for identifying them and rejecting them. 

\medskip

Recall that a coherent sheaf ${\mathcal F}$ on $\PP^n$  is {\bf (Castelnuovo-Mumford) regular} if: 
\[ \mbox{H}^i(\PP^n, {\mathcal F}(-i)) = 0 \ \mbox{for all} \ i > 0 \] 
and its relevance to our discussion is: 

\medskip

\nt {\bf Lemma 5.1.} ${\mathcal F}$ is regular if and only if it is an object of ${\mathcal C}$. \cite{Laz}

\medskip

We use this to detect red dots with the following:

\medskip

\nt {\bf Lemma 5.2.} If ${\mathcal F}$ is a torsion-free sheaf of rank $r$ on $\PP^2$ that is Castelnuovo-Mumford regular and $Z$-semi-stable, then:
\[ \chi({\mathcal F}) \le m_r(\chi'({\mathcal F}))       \] 
where $m_r(n)$ is the value taken by any torsion-free sheaf of rank $r$ and $\chi'({\mathcal F}) = n$ that minimizes the discriminant:
\[\Delta({\mathcal F}) = c_1({\mathcal F}) ^2 -  2\ch_0({\mathcal F}) \ch_2({\mathcal F}) \ge 0  \]

{\bf Proof.} In terms of the Chern classes of ${\mathcal F}$,  
\[ \chi'({\mathcal F}) = \frac 32 \rk({\mathcal F}) + \ch_1({\mathcal F}) \cdot H  \] 
 so if the rank and $\chi'$ are fixed, then $\ch_1({\mathcal F})$ is also fixed, and:
 \[ \chi({\mathcal F}) =   \rk({\mathcal F}) +  \frac 32 \ch_1({\mathcal F}) \cdot H +\ch_2({\mathcal F}) \] 
 shows that $\chi$ is maximized when $\ch_2$ is maximized, i.e. when $\Delta$ is minimized.  
  \qed

\medskip

The following table gives $m_r(\chi')$ for some indicated values of $r$ and $\chi'$ and sets up the next
diagram (with all red dots suppressed). 

\bigskip

\begin{center}

\begin{tabular}{|c|c|c|c|c|c|c|c|c|c|c|} \hline
$r\backslash \chi'$  & $1/2$ & $1$ & $3/2$ & $2$ & $5/2$ & $3$ & $7/2$ & 4 & $9/2$  \\ \hline 
$1$ & $0$ & & $1$ & & $3$ & & $6$ & & $10$  \\ \hline
$2$ & &  & & $0$ & & $2$ & & $3$  & \\ \hline
$3$ & & & & & & & $0$ & & $3$ \\ \hline

\end{tabular} 

\end{center} 

\medskip 

\begin{center}

\begin{tikzpicture} 

\draw[shift={(0,0)},color=black]  node[below] {$(0,0)$};

\draw[->]  (0,0) -- (-1,0);
\draw[shift={(-2,0)},color=black] node[above] {${\mathcal O}(-1)[1]$};
\draw[shift={(-2,0)},color=black] node[below] {${\mathcal O}(-2)[2]$};

\draw[->]  (0,0) -- (0,1);
\draw[shift={(0,1)},color=black] node[above] {$\CC_p$};

\draw[->]  (0,0) -- (3,1);
\draw[shift={(3,1)},color=black] node[right] {${\mathcal O}$};
\draw[shift={(3,1)},color=black] node{$\bullet$};

\draw[->]  (0,0) -- (5,3);
\draw[shift={(5,3)},color=black] node[right] {$\mathcal O(1)$};
\draw[shift={(5,3)},color=black] node{$\bullet$};

\draw[->]  (0,0) -- (5,2);
\draw[shift={(5,2)},color=black] node[left] {$\mathcal I_p(1)$};
\draw[shift={(5,2)},color=black] node{$\bullet$};

\draw[->]  (0,0) -- (6,2);
\draw[shift={(6,2)},color=black] node[right] {${\mathcal O}^2$};
\draw[shift={(6,2)},color=black] node{$\bullet$};

\draw[->]  (0,0) -- (7,6);
\draw[shift={(7,6)},color=black] node[right] {${\mathcal O}(2)$};

\draw[->]  (0,0) -- (7,5);
\draw[shift={(7,5)},color=black] node[right] {$I_p(2)$};
\draw[shift={(7,5)},color=black] node{$\bullet$};

\draw[->]  (0,0) -- (7,4);
\draw[shift={(7,4)},color=black] node[right] {$I_{pq}(2)$};
\draw[shift={(7,4)},color=black] node{$\bullet$};

\draw[->]  (0,0) -- (7,3);
\draw[shift={(7,3)},color=black] node[right] {$I_{pqr}(2)$};
\draw[shift={(7,3)},color=black] node{$\bullet$};

\end{tikzpicture}  

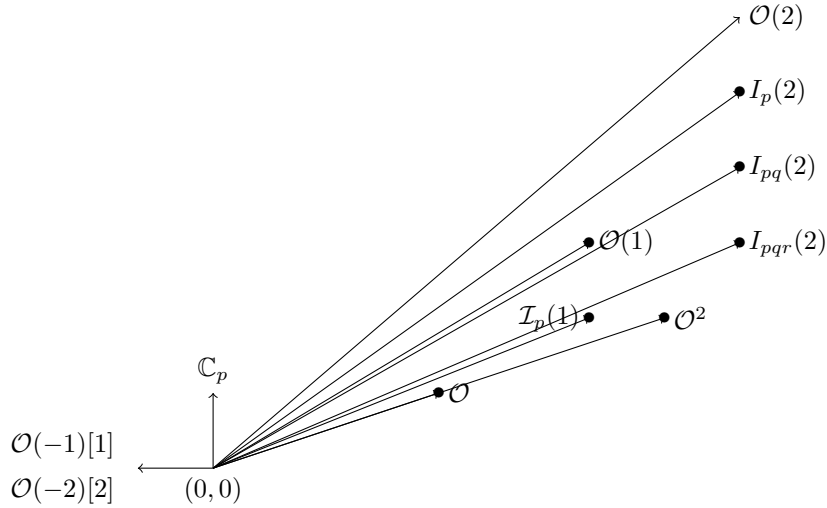
Image of $Z$ for $n = 2$ and $e = 2$ 
\end{center}

Each of the nodes represents a maximal slope destabilizer of an $E(\sigma)$, but some 
additional care needs to be taken to account for filtrations with several terms. 

\medskip

The partition (but not a stratification!) of $\PP^2_4 = \PP^{14}$ is given by open sets of: 

\medskip

$\bullet$ the four-uple embedding of $\PP^2$, which is the closed subset of maps: 
\[ {\mathcal O}(2) \rightarrow {\mathcal O}_p \rightarrow \omega_{\PP^2}(-2)[2] \] 
with Hilbert function $[1 \ 1 \ 1\ 1 \ 1]$. This has dimension 2. 

\medskip

$\bullet$ the variety of secant lines, which has a longer filtration: 
\[ {\mathcal O}_{\PP^2}(2) \rightarrow {\mathcal O}_l(2) \rightarrow {\mathcal O}_{pq} \rightarrow \omega_{\PP^2}(-2)[2] \]
with Hilbert function $[1\ 2\ 2\ 2 \ 1]$. Dimension  $5$.

\medskip

$\bullet$ the variety of linear spans of lines (embedded as rational quartics):
\[ {\mathcal O}_{\PP^2}(2) \rightarrow {\mathcal O}_l(2) \rightarrow \omega_{\PP^2}(-2)[2] \]
and Hilbert function $[1\ 2\ 3\ 2\ 1]$. Dimension $6$.

\medskip

$\bullet$ the variety of linear spans of three non-collinear (in $\PP^2$) points: 
\[ {\mathcal O}_{\PP^2}(2) \rightarrow {\mathcal O}_{pqr} \rightarrow \omega_{\PP^2}(-2)[2] \]
with Hilbert function $[1\ 3\ 3\ 3\ 1]$. Dimension $8$. 

\medskip

$\bullet$ the variety of linear spans of a rational quartic and a point on $\PP^2$:
\[ {\mathcal O}_{\PP^2}(2) \rightarrow {\mathcal O}_{l\cup p}(2) \rightarrow \omega_{\PP^2}(-2)[2] \]
has the same Hilbert function $[1\ 3\ 4\ 3\ 1]$. Dimension $9$.

\medskip

$\bullet$ the variety of linear spans of an intersection of two conics, which is: 
\[ {\mathcal O}_{\PP^2}(2) \rightarrow {\mathcal O}_{\PP^2}(2)/{\mathcal O}_{\PP^2}^2 \rightarrow \omega_{\PP^2}(-2)[2] \]
and is the first case where the complex quotient $Q(\sigma)$  is not a coherent sheaf. 
This also doubles up the Hilbert function  of $[1 \ 3\ 4\ 3\ 1]$. Dimension $11$.

\medskip

$\bullet$ the variety of linear spans of a single conic (embedded as a curve of degree $8$) 
\[ {\mathcal O}_{\PP^2}(2) \rightarrow {\mathcal O}_C(2) \rightarrow  \omega_{\PP^2}(-2)[2] \]
is a hypersurface and has Hilbert function $[1\ 3\ 5\ 3\ 1]$. Dimension $13$.

\medskip

$\bullet$ this leaves the open variety of semi-stable $E(\sigma)$: 
\[ [{\mathcal O}(-2)^7 \rightarrow {\mathcal O}(-1)^7 ] \] 
which are skew-symmetric maps of constant corank one. Dimension $14$. 

\medskip

The two locally closed sets with Hilbert polynomial $[1\ 3\ 4\ 3\ 1]$ are distinguished by 
their betti tables, whose interior squares are, respectively:

\medskip

\begin{center}

\begin{tabular}{|c|c|} \hline
0 & 0 \\ \hline
2 & 1 \\ \hline
 2 & 2 \\ \hline
 1 & 2 \\ \hline
 0 & 0 \\ \hline
 \end{tabular} 
 \ \ \ and \ \ \ 
 \begin{tabular}{|c|c|} \hline
0 & 0 \\ \hline
2 & 0 \\ \hline
 1 & 1 \\ \hline
 0 & 2 \\ \hline
 0 & 0 \\ \hline
 \end{tabular} 
 
\end{center} 

\section{The Odd Degree Case $d = 2e - 1$}

In order to achieve the symmetry in this case, we choose:
\[ Z_{-\frac 12} (E) = ( \chi'_{-\frac 12}(E(t)), \chi_{-\frac 12} (E(t))) \] 
which still has the property that 
$0 < \arg(Z(C) < 1 \ \mbox{when $C$ is a complex in ${\mathcal C}$}$, 
and therefore still defines a stability condition on complexes with:
\[ \chi_{-\frac 12} ({\mathcal O}_{\PP^n}(e)) = \chi_{-\frac 12}({\omega}_{\PP^n}(1-e)[n]) \ \mbox{and} 
\  \chi'_{-\frac 12} ({\mathcal O}_{\PP^n}(e)) = -\chi'_{-\frac 12}({\omega}_{\PP^n}(1-e)[n]) \]

\medskip

For example, for low values of $n$:

\medskip

($n = 1$) 
\[ Z_{-\frac 12}({\mathcal O}_{\PP^1}) = (1, \frac 12), \ Z_{-\frac 12}({\mathcal O}_{\PP^1}(-1)[1]) = (-1, \frac 12) \]  

($n=2$)
\[ Z_{-\frac 12}({\mathcal O}_{\PP^2}) = (1, \frac 38), \ Z_{-\frac 12}({\mathcal O}(-1)_{\PP^2}[1]) = (0, \frac 18), \ Z_{-\frac 12}({\mathcal O}_{\PP^2}(-2)[2]) = (-1,\frac 38) \] 

($n = 3$) 
\[ Z_{-\frac 12}({\mathcal O}_{\PP^3}) = (\frac{23}{24}, \frac{5}{16}), Z_{-\frac 12}({\mathcal O}_{\PP^3}(-1)[1]) = (\frac 1{24}, \frac 1{16}) \]

Note that as before, 
\[ Z_{-\frac 12}(\CC_p) = (0,1) \]

\nt {\bf Proposition 1.} $\CC_p$ and all line bundles on $\PP^n$  are stable for this central charge.  

\medskip

{\bf Proof.} We use more information about the Eagon-Northcott complex:
\[ [EN_\bullet \otimes {\mathcal O}_{\PP^n}(-\bullet) ] = {\mathcal O}_{\PP^n}(d) \] 
namely, that {\bf each} of the representations $EN_i$ is irreducible. It then follows that a maximal destabilizing
subcomplex (if ${\mathcal O}_{\PP^n}(d)$ were unstable) would necessarily be a {\bf truncation}
$ [ EN_{\bullet} \otimes {\mathcal O}_{\PP^n}(-\bullet) ]_{\bullet \le i}$.

\medskip 

Similarly, the complex for $\CC_p$ is a Koszul complex of irreducible representations for the 
subgroup GL$(W) \subset \mbox{GL}(V)$ fixing $p$:
\[ [\wedge^\bullet W \otimes {\mathcal O}_{\PP^n}(-\bullet)] = \CC_p \] 
and any maximal destabilizer is therefore a truncation.

\medskip

As a result, $\CC_p$ and all the line bundles are stable for {\bf any} central charge:
\[ \zeta ({\mathcal O}_{\PP^n}(-i)[i]) := \zeta_i \ \mbox{with} \  \arg(\zeta_n) > \arg(\zeta_{n-1}) > \cdots > \arg(\zeta_0) \] 
since truncations only eliminate vectors with larger arguments than the remainder.

\medskip

The Proposition then follows from the ordering of the (anti)-slopes: 
\[ \nu(t) := \frac{\chi'({\mathcal O}_{\PP^n}(t))}{\chi({\mathcal O}_{\PP^n}(t))} = \sum_{i = 1}^n  \frac{1}{t + i} \] 
\[ \nu(-n-\frac 12) < \nu(-n + \frac 12) < \cdots < \nu(-\frac 12)  \qed \] 

This allows us to proceed as in the $d = 2e$ case, however note that: 

\medskip

(a) {\bf No} socle 
$\sigma: {\mathcal O}_{\PP^n}(e) \rightarrow \omega_{\PP^n}(1-e)[n]$ is injective (dimension count)

\medskip

and, correspondingly, 

\medskip

(b) {\bf Every} cone $E(\sigma)$ over a socle $\sigma$ is unstable. 

\medskip

In fact, since:
\[ \arg(Z_{-\frac 12}({\mathcal O}_{\PP^n}(-i)[i]) < 1 \ \mbox{for all $i$} \] 
it follows that every socle map $\sigma: {\mathcal O}_{\PP^n}(e) \rightarrow \omega_{\PP^n}(1-e) = ({\mathcal O}_{\PP^n}(e))^\vee $ is symmetric:
\[ \ker(\sigma) \rightarrow {\mathcal O}_{\PP^n}(e) \rightarrow Q(e)  = Q(e)^\vee \rightarrow ({\mathcal O}_{\PP^n}(e))^\vee  \rightarrow \ker(\sigma)^\vee \] 
(and the identification of $Q(e)$ with $Q(e)^\vee$ is skew symmetric). 

\medskip

There are no surprises in the $n = 1$ case. The values merely shift:

\begin{center}

\begin{tikzpicture} 

\draw[shift={(0,0)},color=black]  node[below] {$(0,0)$};

\draw[->]  (0,0) -- (-1,.5);
\draw[shift={(-1,.5)},color=black] node[left] {${\mathcal O}(-1)[1]$};

\draw[->]  (0,0) -- (0,1);
\draw[shift={(0,1)},color=black] node[above] {$\CC_p$};

\draw[->]  (0,0) -- (2,0);
\draw[shift={(2,0)},color=black] node[right] {$E(\sigma)$};

\draw[->]  (0,0) -- (1,.5);
\draw[shift={(1,.5)},color=black] node[right] {$\mathcal O$};
\draw[shift={(1,.5)},color=black] node{$\bullet$};

\draw[->]  (0,0) -- (1,1.5);
\draw[shift={(1,1.5)},color=black] node[right] {$\mathcal O(1)$};
\draw[shift={(1,1.5)},color=black] node{$\bullet$};

\draw[->]  (0,0) -- (1,3.5);
\draw[shift={(1,3.5)},color=black] node[right] {$\mathcal O(e)$};

\draw[shift={(1,2.5)},color=black] node[right] {$\vdots$};
\end{tikzpicture}  

Image of $Z_{-\frac 12}$ for $n = 1$ and $d = 2e - 1$ 

\end{center} 

\medskip

The arrows and nodes for $n = 2$ and $d = 1$ are: 

\begin{center}

\begin{tikzpicture}

\draw[->]  (0,0) -- (-1,.75);
\draw[shift={(-1,.75)},color=black] node[left] {${\mathcal O}(-2)[2]$};

\draw[->]  (0,0) -- (0,.25);
%\draw[shift={(0,.25)},color=black] node{$\bullet$};

\draw[->]  (0,0) -- (0,2);
\draw[shift={(0,2)},color=black] node[above]{$\CC_p$};

\draw[->]  (0,0) -- (1,.75);
\draw[shift={(1,.75)},color=black] node[below] {${\mathcal O}$};
\draw[shift={(1,.75)},color=red] node{$\bullet$};

\draw[->]  (0,0) -- (2,1.5);
\draw[shift={(2,1.5)},color=black] node[below] {${\mathcal O}^2$};
\draw[shift={(2,1.5)},color=red] node{$\bullet$};

\draw[->]  (0,0) -- (2,3.75);
\draw[shift={(2,3.75)},color=black] node[right] {${\mathcal O}(1)$};

\draw[->]  (0,0) -- (2,1.75);
\draw[shift={(2,1.75)},color=black] node[right] {$I_p(1)$};
\draw[shift={(2,1.75)},color=black] node {$\bullet$};

\end{tikzpicture}  

Image of $Z_{-\frac 12}$ for $n = 2$ and $e = 1$ 
\end{center}

One novelty of this odd case is the arrow:
\[ Z_{-\frac 12}({\mathcal O}_{\PP^2}(-1)[1]) = (0,\frac 18) \] 
and another is the observation that: 

\medskip

$\bullet$ for all $e$, the kernel ker$(\sigma)$ contains $\CC^{e+1} \otimes {\mathcal O}_{\PP^2}$

\medskip

This is responsible for the red dot at $Z_{-\frac 12}({\mathcal O})$ but also at $Z_{-\frac 12}({\mathcal O}^2)$ since:

\medskip

$\bullet$ every map ${\mathcal O}^2 \rightarrow {\mathcal O}(1)$ {\bf factors through} ${\mathcal I}_p(1)$ (with a kernel!)

\medskip

And finally, the first interesting case $d = 3$ with all red dots but one surpressed:

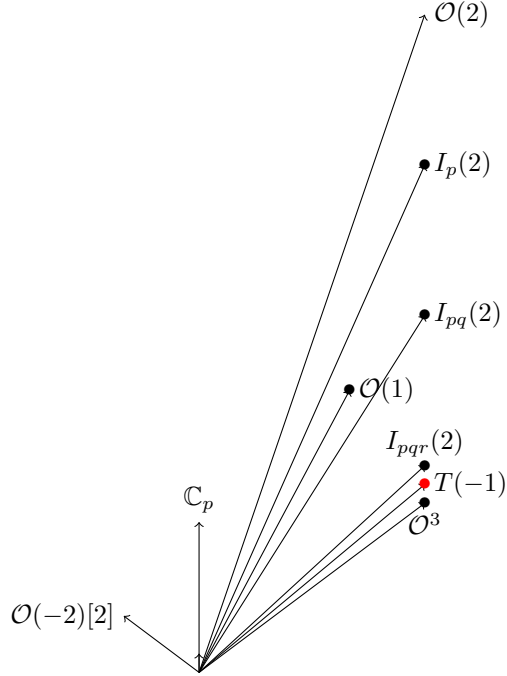
\begin{center}

\begin{tikzpicture}

\draw[->]  (0,0) -- (-1,.75);
\draw[shift={(-1,.75)},color=black] node[left] {${\mathcal O}(-2)[2]$};

\draw[->]  (0,0) -- (0,.25);
%\draw[shift={(0,.25)},color=black] node{$\bullet$};

\draw[->]  (0,0) -- (0,2);
\draw[shift={(0,2)},color=black] node[above]{$\CC_p$};

\draw[->]  (0,0) -- (2,3.75);
\draw[shift={(2,3.75)},color=black] node[right] {${\mathcal O}(1)$};
\draw[shift={(2,3.75)},color=black] node{$\bullet$};

\draw[->]  (0,0) -- (3,8.75);
\draw[shift={(3,8.75)},color=black] node[right] {${\mathcal O}(2)$};

\draw[->]  (0,0) -- (3,6.75);
\draw[shift={(3,6.75)},color=black] node[right] {$I_p(2)$};
\draw[shift={(3,6.75)},color=black] node {$\bullet$};

\draw[->]  (0,0) -- (3,4.75);
\draw[shift={(3,4.75)},color=black] node[right] {$I_{pq}(2)$};
\draw[shift={(3,4.75)},color=black] node {$\bullet$};

\draw[->]  (0,0) -- (3,2.5);
\draw[shift={(3,2.5)},color=black] node[right] {$T(-1)$};
\draw[shift={(3,2.5)},color=red] node {$\bullet$};

\draw[->]  (0,0) -- (3,2.75);
\draw[shift={(3,2.75)},color=black] node[above] {$I_{pqr}(2)$};
\draw[shift={(3,2.75)},color=black] node {$\bullet$};

\draw[->]  (0,0) -- (3,2.25);
\draw[shift={(3,2.25)},color=black] node[below] {${\mathcal O}^3$};
\draw[shift={(3,2.25)},color=black] node{$\bullet$};

\end{tikzpicture}  

Image of $Z_{-\frac 12}$ for $n = 2$ and $d = 3$ 

\end{center}

As with $d = 1$, the red node arises because each map from the tangent bundle:
\[ T_{\PP^2}(-1) \rightarrow {\mathcal O}_{\PP^2}(2) \ \mbox{factors through} \ {\mathcal I}_{pqr}(2) \] 
due to the fact that the the second chern class $c_2(T) = 3$. This addresses Example 3.2 in a novel way. 
and in this case, 
\[ {\mathcal O}_{\PP^2}^3 \rightarrow {\mathcal O}_{\PP^2}(2) \ \mbox{need not factor through any intermediary} \] 
and therefore accounts for the open stratum with betti table: 
 
\[ \begin{tabular}{|c|c|c|c|c|} \hline
 1 & 0 & 0 & 0   \\ \hline
 0 & 3 & $0$ & 0  \\ \hline
 0 & $0$ & 3 &  0  \\ \hline  
 0 & 0 & 0 & 1 \\ \hline
 \end{tabular} \] 

\nt while the kernel $I_{pqr}(2)$ (for non-collinear points) accounts for: 

\[ \begin{tabular}{|c|c|c|c|c|} \hline
 1 & 0 & 0 & 0   \\ \hline
 0 & 3 & $2$ & 0  \\ \hline
 0 & $2$ & 3 &  0  \\ \hline  
 0 & 0 & 0 & 1 \\ \hline
 \end{tabular} \] 

\nt It is now an exercise to match the rest of the betti tables with this partition. 

\bigskip

\nt {\bf Final Remark.} This may be useful for higher values of $n$.  The first test is to use it on Example 3.3. As a starting point, one needs to 
thoroughly understand the ``destabilizing'' 
subcomplexes of an Eagon-Northcott complex without the crutch (in $n = 2$) of knowing they are coherent sheaves so that arrows and nodes can be compiled. A good starting point for 
this could be the recent result \cite{BS}.

\end{document}